\documentclass[12pt]{article}
\usepackage{amsthm,amsfonts,amsmath,amscd,amssymb,latexsym}    
\usepackage{mathrsfs}
\usepackage{epsfig,graphicx,color}     
\usepackage[all]{xy}  
\usepackage{makeidx}
\newif\ifprintable
\printabletrue
\ifprintable

\newcommand{\url}[1]{{#1}}

\else
\usepackage[colorlinks=true]{hyperref}
\fi
\makeindex
\newif\ifmarking
\markingtrue

\title{Uniqueness of Symplectic Structures}

\author{
Dietmar Salamon\thanks{Partially supported by the 
Swiss National Science Foundation Grant 200021-127136}\\
ETH Z\"urich
}

\date{12 December 2012}

\newtheorem{PARA}{}[section]

\newtheorem{example}[PARA]{Example} 
\newtheorem{conjecture}[PARA]{Conjecture}

\newcommand{\MAT}[1]{\left[\begin{array}{#1}}
\newcommand{\RIX}{\end{array}\right]}
\newcommand{\p}{\partial}
\newcommand{\one}{{{\mathchoice {\rm 1\mskip-4mu l} {\rm 1\mskip-4mu l}
{\rm 1\mskip-4.5mu l} {\rm 1\mskip-5mu l}}}}

\newcommand{\C}{{\mathbb C}}

\newcommand{\R}{{\mathbb R}}
\newcommand{\T}{{\mathbb T}}
\newcommand{\Z}{{\mathbb Z}}

\newcommand{\cC}{{\mathcal C}} 

\newcommand{\cE}{{\mathcal E}}

\newcommand{\cJ}{{\mathcal J}}
\newcommand{\cK}{{\mathcal K}}

\newcommand{\cR}{{\mathcal R}}
\newcommand{\cS}{{\mathcal S}}

\newcommand{\sF}{\mathscr{F}}

\newcommand{\Om}{{\Omega}}
\newcommand{\om}{{\omega}}
\newcommand{\eps}{{\varepsilon}}
\renewcommand{\phi}{{\varphi}}

\newcommand{\im}{{\mathrm{im}}}

\newcommand{\id}{{\mathrm{id}}}

\newcommand{\standard}{{\mathrm{std}}}

\newcommand{\Aut}{{\mathrm{Aut}}}   
 
\newcommand{\Diff}{{\mathrm{Diff}}}

\newcommand{\FS}{{\mathrm{FS}}}
\newcommand{\SW}{{\mathrm{SW}}}

\newcommand{\dvol}{{\rm dvol}}

\newcommand{\symp}{{\mathrm{symp}}}

\newcommand{\SO}{{\mathrm{SO}}} 
\newcommand{\Spin}{{\mathrm{Spin}}}

\newcommand{\PU}{{\mathrm{PU}}}

\renewcommand{\o}{{\frak o}}

\newcommand{\CP}{{\C\mathrm{P}}}

\newcommand{\HF}{{\mathrm{HF}}}

\newcommand{\PD}{{\mathrm{PD}}}

\newcommand{\inner}[2]{\langle #1, #2\rangle}

\def\NABLA#1{{\mathop{\nabla\kern-.5ex\lower1ex\hbox{$#1$}}}}
\def\Nabla#1{\nabla\kern-.5ex{}_{#1}}
\def\abs#1{\mathopen|#1\mathclose|}
\def\Abs#1{\left|#1\right|}

\def\Norm#1{\left\|#1\right\|}

\renewcommand{\p}{{\partial}}

\newcommand{\longhookrightarrow}{\ensuremath{\lhook\joinrel\relbar\joinrel\rightarrow}}
\newcommand{\notimplies}{{\hspace{7pt}\not\hspace{-7pt}\implies}}

\begin{document}
\maketitle

\begin{abstract}
This survey paper discusses some uniqueness 
questions for symplectic forms on compact 
manifolds without boundary.
\end{abstract}


\section{Introduction}

Let $M$ be an oriented $2n$-manifold. A symplectic form on $M$ is a 
closed $2$-form $\om\in\Om^2(M)$ whose top exterior power $\om^n$ 
is a volume form.  Fundamental questions in symplectic topology are 
which deRham cohomology classes $a\in H^2(M;\R)$ are represented
by symplectic forms (the existence problem) and whether any two
symplectic forms on $M$ that represent the same cohomology class
are related by a suitable equivalence relation (the uniqueness problem).
Necessary conditions for the existence of symplectic forms are the
existence of an almost complex structure and, in the closed case, 
the existence of a cohomology class $a\in H^2(M;\R)$ with $a^n>0$.
Whether or not these conditions are also sufficient, when $M$ 
is a closed manifold, is completely open in dimensions $2n\ge 6$.
In dimension four additional necessary conditions for existence arise 
from Seiberg--Witten theory.  In the present survey paper the focus is mainly 
on the uniqueness problem.  Section~\ref{sec:SYMP} gives precise 
formulations of some relevant questions, Section~\ref{sec:EX}
discusses what is known about these question in some examples
(without any claim of providing a complete picture), and Section~\ref{sec:D}
discusses some conclusions and conjectures as well as the Donaldson geometric 
flow approach to the uniqueness problem for hyperk\"ahler surfaces.


\section{The space of symplectic forms}\label{sec:SYMP}

\begin{PARA}[{\bf Equivalence relations}]
\label{para:equiv}\rm
Let $M$ be a manifold without boundary.
Consider the following statements for two symplectic 
forms $\om_0$, $\om_1$ on $M$.
\begin{description}
\item[(a)] 
$\om_0$ and $\om_1$ are connected 
by a path of cohomologous symplectic forms.
\item[(b)]
$\om_0$ and $\om_1$ are connected 
by a path of symplectic forms.
\item[(c)]
$\om_0$ and $\om_1$ are connected 
by a path of nondegenerate $2$-forms.
\item[(d)]
$\om_0$ and $\om_1$ have the same first Chern class in $H^2(M;\Z)$.
\item[(A)]
There exists a diffeomorphism $\phi$ of $M$ such that
$\om_0=\phi^*\om_1$.
\item[(B)]
There exists a diffeomorphism $\phi$ of $M$ such that
$\om_0$ and $\phi^*\om_1$ are connected 
by a path of symplectic forms.
\item[(C)]
There exists a diffeomorphism $\phi$ of $M$ such that
$\om_0$ and $\phi^*\om_1$ are connected 
by a path of nondegenerate $2$-forms.
\item[(D)]
There exists a diffeomorphism $\phi$ of $M$ such that
$c_1(\om_0)=\phi^*c_1(\om_1)$.
\end{description}
\bigbreak

\noindent
These statements define equivalence relations.
Two symplectic forms~$\om_0$ and~$\om_1$ are called 
{\bf isotopic} if they are related by~(a) 
(i.e.\ a path of cohomologous symplectic forms),
they are called {\bf homotopic} if they are related by~(b)
(i.e.\ a path of symplectic forms),
they are called {\bf diffeomorphic} 
if they are related by~(A) (i.e.\ a diffeomorphism),
and they are called {\bf deformation equivalent} 
if they are related by~(B) (i.e.\ a diffeomorphism, 
followed by a path of symplectic forms). 

For closed manifolds these equivalence relations 
are related as follows
$$
\begin{array}{ccccccccc}
(a) & \implies & (b) & \implies & (c) & \implies & (d) \\
\Downarrow & & \Downarrow & & \Downarrow & & \Downarrow \\
(A) & \implies & (B) & \implies & (C) & \implies & (D)
\end{array}.
$$
In particular, by Moser isotopy, assertion~(a) holds 
if and only if there exists a diffeomorphism $\phi:M\to M$ 
isotopic to the identity such that $\phi^*\om_1=\om_0$.
Thus~(a) implies~(A).  For~(A), (B), (C) it may also be interesting 
to restrict to diffeomorphisms that act as the identity on homology
(see Karshon--Kessler--Pinsonnault~\cite{KKP}).
For~(A) this would mean that equivalent symplectic forms represent 
the same cohomology class and for~(B) that they have the same first 
Chern class and the same Gromov--Witten invariants.
\end{PARA}

\begin{PARA}[{\bf Three questions}]\label{para:questions}\rm
The uniqueness problem in symplectic topology is the problem 
of understanding these equivalence relations on the space of 
symplectic forms. Here are some questions that may serve  
as a guideline. Fix a deRham homology class 
$$
a\in H^2(M;\R)
$$
and denote the space of symplectic forms 
on $M$ representing the class $a$ by
$$
\cS_a := \left\{\rho\in\Om^2(M)\,|\,
\rho^n\ne 0,\,d\rho=0,\,[\rho]=a\right\}.
$$
This is an open set in the Fr\'echet space of all
closed $2$-forms on $M$ representing the class $a$.
The following three questions correspond to the equivalence 
relations~(a), (A), (B).  In all three questions assume that 
$\cS_a$ is nonempty (respectively that $M$ admits a symplectic form).

\medskip\noindent{\bf Question~1.}
{\it Is $\cS_a$ connected?}

\medskip\noindent{\bf Question~2.}
{\it Are any two symplectic forms in $\cS_a$ diffeomorphic?}

\medskip\noindent{\bf Question~3.}
{\it Are any two symplectic forms on $M$
deformation equivalent?}

\medskip\noindent
Section~\ref{sec:EX} discusses what is known about these questions 
in some examples.  Section~\ref{sec:D} includes a brief discussion 
of the Donaldson geometric flow approach to Question~1 
in dimension four (see~\cite{D4}).
\end{PARA}

\begin{PARA}[{\bf Almost complex structures}]\label{para:ac}\rm
Denote the set of almost complex structures
on $M$ that are compatible with a nondegenerate
$2$-form $\rho$ by 
$$
\cJ(M,\rho) := \left\{J:TM\to TM\,|\,J^2=-\one,\,
\rho(\cdot,J\cdot)\mbox{ is a Riemannian metric}\right\}.
$$
This space is nonempty and contractible.  
In fact the space of nondegenerate $2$-forms on $M$
is homotopy equivalent to the space of almost complex 
structures on $M$. (The homotopy equivalence
is a smooth map $\rho\mapsto J_\rho$, depending 
on the choice of a background metric, which assigns to 
every nondegenerate $2$-form $\rho$ an almost complex 
structure $J_\rho\in\cJ(M,\rho)$; see~\cite{GROMOV2,MS1}.)
Thus every nondegenerate $2$-form $\rho$ on $M$ 
determines a cohomology class 
$$
c_1(\rho):=c_1(TM,J)\in H^2(M;\Z),\qquad J\in\cJ(M,\rho),
$$
called the {\bf first Chern class of $\rho$}. Consider the set
$$
\cC_a := \left\{c_1(\rho)\,|\,\rho\in\cS_a\right\}
\subset H^2(M;\Z).
$$
If this set contains more than one element, 
the answer to Question~1 is negative. 
Moreover, define
\begin{equation*}
\begin{split}
\cC_\symp 
&:= \left\{c_1(\rho)\,|\,\rho
\mbox{ is a symplectic form on }M\right\}
= \bigcup_a\cC_a, \\
\cC 
&:= \left\{c_1(\rho)\,|\,\rho
\mbox{ is a nondegenerate $2$-form on }M\right\}.
\end{split}
\end{equation*}

\smallskip\noindent
{\bf Further questions:}
What is the set $\cC_a$?  Is $\cC_\symp=\cC$?
Does the diffeomorphism group act transitively on $\cC_\symp$?
\end{PARA}

\begin{PARA}[{\bf Homotopy classes of almost complex structures}]
\label{para:homotopyJ}\rm
Let $M$ be a compact connected oriented smooth four-manifold
without boundary.  Denote by $\chi$ the Euler characteristic 
and by $\sigma$ the signature of~$M$.  
By a theorem of Wu an integral cohomology class $c\in H^2(M;\Z)$ 
is the first Chern class of an almost complex structure on $M$,
compatible with the orientation, if and only if 
it is an integral lift of the second Stiefel--Whitney class and 
\begin{equation}\label{eq:hirzebruch}
c^2=2\chi+3\sigma.  
\end{equation}
(The necessity of equation~\eqref{eq:hirzebruch} 
is the Hirzebruch signature theorem.)
Given such a class $c$, denote by
$$
\cJ(M,c) := \left\{J\in\Aut(TM)\,\bigg|\,
\begin{array}{l}
J^2=-\one,\,c_1(TM,J)=c,\\
J\mbox{ induces the given orientation}
\end{array}
\right\}
$$ 
the space of almost complex structures on $M$ with first Chern 
class $c$ and compatible with the orientation of $M$.  
If~$M$ is simply connected then $\cJ(M,c)$ has precisely 
two connected components for every $c\in\cC$.  In general, 
there is a bijection
\begin{equation}\label{eq:mrowka}
\pi_0(\cJ(M,c))\cong 
\mathrm{Tor}_2(H^2(M;\Z))\times
\left(\Z/2\Z\oplus\frac{H^3(M;\Z)}{H^1(M;\Z)\cup c}\right)
\end{equation}
for every $c\in\cC$. (This was pointed out to me a long time ago by
Tom Mrowka, but I don't know a reference.)
The first factor is the set
$$
\mathrm{Tor}_2(H^2(M;\Z)):=\left\{a\in H^2(M;\Z)\,|\,2a=0\right\}.
$$
It characterizes the isomorphism classes of spin$^c$ 
structures on $TM$ with first Chern class $c$ 
(see~\ref{para:SW} and Lawson--Michelsohn~\cite[Appendix~D]{LM}).  
The second factor characterizes the set of homotopy 
classes of almost complex structures on $M$ 
whose canonical spin$^c$ structures are isomorphic to 
a given spin$^c$ structure with first Chern class $c$.  
(This can be proved with a standard Pontryagin manifold 
type construction as in Milnor~\cite[\S 7]{MILNOR}.)
\end{PARA}

\begin{PARA}[{\bf The existence problem}]\label{para:existence}\rm
A natural question is under which conditions 
a given cohomology class $a\in H^2(M;\R)$ can be represented
by a symplectic form. This is the fundamental existence problem 
in symplectic topology. The obvious necessary conditions 
are the existence of an almost complex structure on $M$ and, 
in the case of a closed manifold, the condition $a^n\ne0$. 
In the open case a theorem of Gromov~\cite{GROMOV0,GROMOV2} asserts 
that the existence of an almost complex structure is also sufficient.
In the closed case there are counterexamples 
in dimension four, based on Seiberg--Witten theory
(see Taubes~\cite{TAUBES2}).  The simplest example 
is $M=\CP^2\#\CP^2\#\CP^2$ which does not admit 
a symplectic form because its Seiberg--Witten invariants vanish.  
In higher dimensions the existence problem is completely open.
\end{PARA}

\begin{PARA}[{\bf Diffeomorphism groups}]\label{para:diffeo}\rm
It is also interesting to investigate the topology of the space 
$\cS_a$ of symplectic forms in a given cohomology class 
beyond the question of connectedness.  This is closely related 
to the topo\-logy of certain diffeomorphism groups.
Let $(M,\om)$ be a compact symplectic manifold without 
boundary and denote $a:=[\om]\in H^2(M;\R)$.
Consider the diffeomorphism groups
\begin{equation*}
\begin{split}
\Diff_0(M) &:= 
\left\{\phi\in\Diff(M)\,|\,
\phi\mbox{ is isotopic to the identity}\right\}, \\
\Diff_0(M,\om) &:= 
\left\{\phi\in\Diff_0(M)\,|\,
\phi^*\om=\om\right\}.
\end{split}
\end{equation*}
Thus $\Diff_0(M,\om)$ denotes the group of symplectomorphisms
of $M$ that are {\it smoothly} isotopic to the identity.  
Care must be taken to distinguish it from the subgroup of
all symplectomorphisms of $M$ that are {\it symplectically} 
isotopic to the identity (i.e.\ the identity 
component of $\Diff(M,\om)$). Denote by
$$
\cS_\om := \bigl\{\rho\in\cS_a\,|\,
\rho\stackrel{(\mbox{\tiny{a}})}{\sim}\om\bigr\}
$$
the connected component of $\om$ in $\cS_a$.
(Here $\stackrel{(\mbox{\tiny{a}})}{\sim}$ denotes the 
equivalence relation~(a) in~\ref{para:equiv}.)
By Moser isotopy the group $\Diff_0(M)$ 
acts transitively on~$\cS_\om$.
Hence the map 
$
\Diff_0(M)\to\cS_\om:\phi\mapsto\phi^*\om
$
induces a homeomorphism
$$
\Diff_0(M)/\Diff_0(M,\om)\cong \cS_\om.
$$
Thus the topology of $\cS_\om$ is closely 
related to the topology of the symplectomorphism group.
For example, if there exists a symplectomorphism of $(M,\om)$
that is smoothly, but not symplectically, 
isotopic to the identity then $\Diff_0(M,\om)$ 
is not connected and hence $\cS_\om$ has a 
nontrivial fundamental group
(see Example~\ref{ex:seidel}).  
\end{PARA}

\begin{PARA}[{\bf Seiberg--Witten theory}]\label{para:SW}\rm 
Many results about symplectic structures
in dimension four rely on Taubes--Seiberg--Witten theory.
Here is some relevant background.
Let $M$ be a closed oriented smooth four-manifold.

\smallskip\noindent{\bf 1.} 
Denote the set of equivalence classes of spin$^c$ 
structures on $M$ by $\Spin^c(M)$.  
The tensor product with complex line bundles
defines a free and transitive action 
$H^2(M;\Z)\times\Spin^c(M)\to\Spin^c(M):
(e,\Gamma)\mapsto\Gamma_e$. Each $\Gamma\in\Spin^c(M)$ 
has a first Chern class $c_1(\Gamma)\in H^2(M;\Z)$ 
and $c_1(\Gamma_e)=c_1(\Gamma)+2e$.

\smallskip\noindent{\bf 2.} 
Every nondegenerate $2$-form $\rho\in\Om^2(M)$, 
compatible with the orientation, determines a
canonical equivalence class of spin$^c$ structures 
$\Gamma_\rho\in\Spin^c(M)$, associated to the 
$\rho$-compatible almost complex structures.  
It depends only on the homotopy class of $\rho$ 
and has the first Chern class 
$c_1(\Gamma_\rho)=c_1(\rho)$.  

\smallskip\noindent{\bf 3.} 
Every nondegenerate $2$-form $\rho\in\Om^2(M)$,
compatible with the orientation, determines a
{\it homological orientation of $M$}, i.e.\ 
an orientation $\o_\rho$ of the real vector 
space $H^0(M;\R)\oplus H^1(M;\R)\oplus H^{2,+}(M;\R)$ 
(see Donaldson~\cite{D1}).  

\smallskip\noindent{\bf 4.} 
If $b^+=\dim H^{2,+}(M;\R)=1$ 
then the {\em positive cone} 
$$
\cK:=\left\{a\in H^2(M;\R)\,|\,a^2>0\right\}
$$
has two connected components.  If $\rho$ 
is a symplectic form, compatible with the orientation,
let $\kappa_\rho\subset\cK$ be the connected 
component containing $[\rho]$. 

\smallskip\noindent{\bf 5.} 
The Seiberg--Witten invariants of $M$ take the form of
a map 
$$
\Spin^c(M)\to\Z:\Gamma\mapsto\SW(M,\o,\Gamma).
$$
This map depends on the choice of 
a homological orientation $\o$ of $M$ (see~3). 
Changing the homological orientation 
reverses the sign of the invariant.
When $b^+=1$ the Seiberg--Witten invariant 
also depends on the choice of a connected 
component $\kappa$ of $\cK$ (see~4).
Changing the connected component is governed
by the wall crossing formula of Li--Liu~\cite{LL}.

\smallskip\noindent{\bf 6.} 
Let $\rho$ be a symplectic form,
compatible with the orientation, and $a:=[\rho]$. 
When $b^+\ge2$, Taubes proved in~\cite{TAUBES1,TAUBES2} that
\begin{equation}\label{eq:SWT1}
\SW(M,\o_\rho,\Gamma_\rho)=1,
\end{equation}
\begin{equation}\label{eq:SWT2}
\SW(M,\o_\rho,\Gamma_{\rho,e})\ne0
\qquad\implies\qquad a\cdot e\ge 0,
\end{equation}
\begin{equation}\label{eq:SWT3}
\SW(M,\o_\rho,\Gamma_{\rho,e})\ne0
\quad\mbox{ and }\quad a\cdot e=0
\qquad\implies\qquad e=0
\end{equation}
for all $e\in H^2(M;\Z)$.
In~\cite{TAUBES3} Taubes proved that every 
cohomology class~$e$ with $\SW(M,\o_\rho,\Gamma_{\rho,e})\ne0$ 
is Poincar\'e dual to a $\rho$-symplectic submanifold of~$M$.
These results continue to hold when $b^+=1$,
with $\SW(M,\o_\rho,\Gamma_{\rho,e})$ 
replaced by $\SW(M,\o_\rho,\kappa_\rho,\Gamma_{\rho,e})$
(see~\cite{TAUBES1,TAUBES2,TAUBES3,TAUBES4,TAUBES5}).
\end{PARA}

\bigbreak

An immediate corollary of Taubes' theorems 
is the uniqueness of the first Chern class 
for cohomologous symplectic forms in dimension four.

\medskip\noindent{\bf Corollary~A.}
{\it Let $M$ be a closed smooth four-manifold.
Two cohomologous symplectic forms on $M$ 
have equivalent spin$^c$ structures and hence 
have the same first Chern class.}

\begin{proof}
Let $\rho,\rho'$ be cohomologous symplectic 
forms on $M$ and $a:=[\rho]=[\rho']$.
Choose the orientation such that $a^2>0$ and
choose $e$ such that $\Gamma_{\rho'}=\Gamma_{\rho,e}$.
Assume first that $b^+\ge2$. By~\eqref{eq:SWT1}, 
$\SW(M,\o_{\rho'},\Gamma_{\rho,e}) = \SW(M,\o_{\rho'},\Gamma_{\rho'})=1$
and so $\SW(M,\o_\rho,\Gamma_{\rho,e})\ne0$.
Hence $a\cdot e\ge 0$ by~\eqref{eq:SWT2}. Interchanging 
the roles of $\rho$ and $\rho'$ gives $a\cdot e\le 0$,
hence $a\cdot e=0$, and hence $e=0$ by~\eqref{eq:SWT3}.
This shows that $\Gamma_{\rho'}=\Gamma_{\rho,e}=\Gamma_\rho$
and $c_1(\rho')=c_1(\rho)+2e=c_1(\rho)$.  The proof for
$b^+=1$ is verbatim the same.  Include
$\kappa_\rho$, respectively $\kappa_{\rho'}$,
as an argument of $\SW$ and use the fact that 
$\kappa_\rho=\kappa_{\rho'}$, because $[\rho]=[\rho']$.
\end{proof}

The next corollary is a special case of a result by 
Conolly-L\'e-Ono~\cite{CLO}.  It strengthens 
Corollary~A in the simply connected case.

\medskip\noindent{\bf Corollary~B.}
{\it Let $M$ be a simply connected 
closed smooth four-manifold.}

\smallskip\noindent{\bf (i)}
{\it Two symplectic forms on $M$ with the same first 
Chern class and inducing the same orientation on $M$
are homotopic as nondegenerate $2$-forms.} 

\smallskip\noindent{\bf (ii)}
{\it Two cohomologous symplectic forms on $M$ 
are homotopic as nondegenerate $2$-forms.} 

\begin{proof} (Following Conolly--L\'e--Ono~\cite{CLO}.)
Assume $b^+=1$. (The case $b^+\ge2$ is easier.)
Fix an orientation of~$M$ and an integral lift $c\in H^2(M;\Z)$ 
of the second Stiefel--Whitney class with $c^2=2\chi+3\sigma$.
Let $\cR^c$ be the set of nondegenerate $2$-forms on $M$ 
with first Chern class $c$ and compatible 
with the orientation. Then $\cR^c$ has precisely 
two connected components, by~\eqref{eq:mrowka}.
Let $\rho,\rho'\in\cR^c$ be symplectic forms.
Then $c_1(\rho')=c_1(\rho)=c$ and so
$\Gamma_{\rho'}=\Gamma_\rho$, 
since $M$ is simply connected.
By~\eqref{eq:SWT1},
$
\SW(M,\o_\rho,\kappa_\rho,\Gamma_\rho) = 1
= \SW(M,\o_{\rho'},\kappa_{\rho'},\Gamma_{\rho'}).
$
The wall crossing formula in Li--Liu~\cite{LL} asserts, 
in the simply connected case, that 
$\SW(M,\o,\kappa,\Gamma)-\SW(M,\o,-\kappa,\Gamma)=\pm1$
whenever $c_1(\Gamma)^2\ge2\chi+3\sigma$.
Since $\SW(M,\o_\rho,\kappa_\rho,\Gamma_\rho)
-\SW(M,\o_\rho,\kappa_{\rho'},\Gamma_\rho)$ 
is even, it follows that ${\kappa_{\rho'}=\kappa_\rho}$. 
Hence $\SW(M,\o_{\rho'},\kappa_\rho,\Gamma_\rho) 
= 1 = \SW(M,\o_\rho,\kappa_\rho,\Gamma_\rho)$,
and hence $\o_{\rho'} = \o_\rho$.
In~\cite[Prop.~3.25]{D1} and~\cite[Lemma~6.4]{D2} 
Donaldson proved that there is a free involution 
on $\pi_0(\cR^c)$, which reverses the homological 
orientation. Hence any two symplectic forms 
in $\cR^c$ belong to the same connected component.
This proves~(i).  Assertion~(ii) follows 
from~(i) and Corollary~A.
\end{proof}


\section{Examples}\label{sec:EX}

\begin{example}[{\bf Open manifolds}]\label{ex:open}\rm
Let $M$ be a connected noncompact smooth manifold
that admits an almost complex structure.  
Then the h-principle rules.  Namely, a theorem 
of Gromov asserts that, for every de\-Rham cohomology 
class $a\in H^2(M;\R)$, the inclusion of the space 
$\cS_a$ of all symplectic forms representing 
the class $a$ into the space of all nondegenerate $2$-forms 
is a homotopy equivalence. (See~\cite{GROMOV0},
\cite[page~84]{GROMOV2}, \cite[Theorem~10.2.2]{EM}, 
\cite[Theorem~7.34]{MS1}.)  This implies the existence 
statement in~\ref{para:existence} and shows that the
uniqueness problem reduces to topological 
obstruction theory. 

For Euclidean space $M=\R^{2n}$ 
the answer to question~3 is positive and 
the answers to questions~1 and~2 are negative.
Namely, there are two homotopy classes of symplectic 
forms on $\R^{2n}$, one for each orientation, and
in~\cite{GROMOV1} Gromov constructed a symplectic 
form $\om$ on $\R^{2n}$ (for $n\ge 2$) such that 
$(\R^{2n},\om)$ is not symplectomorphic
to any open subset of $(\R^{2n},\om_0)$ with the 
standard symplectic form (see also~\cite{BP} 
and~\cite[Example~13.8]{MS1}).
By contrast, if $\om$ is a symplectic form 
on $\R^4$ that agrees with the standard symplectic 
form $\om_0$ at infinity, then $(\R^4,\om)$ is 
symplectomorphic to $(\R^4,\om_0)$ by another theorem
of Gromov (see~\cite{GROMOV1} and~\cite[Theorem~9.4.2]{MS2}).
\end{example}

\begin{example}[{\bf Closed two-manifolds}]\label{ex:Riemann}\rm
Let $M$ be a closed orientable two-manifold.
Then the space $\cS_a$ is nonempty and convex for every
nonzero cohomology class $a\in H^2(M;\R)$. Since $M$ 
admits an orientation reversing diffeomorphism,
questions~1, 2, and~3 have positive answers.
\end{example}

In higher dimensions the uniqueness problem for symplectic forms 
on closed manifolds does not reduce to topological obstruction 
theory.  There is often a dramatic difference between 
the space of nondegenerate two-forms and the space 
of symplectic forms, as the following examples show.

\begin{example}[{\bf A six-manifold}]\label{ex:McDuff}\rm
This example is due to McDuff~\cite{MCDUFF3}
(see also~\cite[Theorem~9.7.4]{MS2}).
Here the set $\cS_a$ is disconnected.  
However, by an arbitrarily small perturbation of the 
cohomology class $a$ the two known distinct 
connected components of $\cS_a$ merge to a single
connected component.  McDuff's 
paper~\cite{MCDUFF3} also contains a variant of 
this construction in dimension eight. 
Consider the manifold
$
M := \T^2\times S^2\times S^2.
$
Identify the $2$-torus  
with the product of two circles.
For $\theta\in S^1$ and $z\in S^2$ let
$\phi_{z,\theta}:S^2\to S^2$ be the rotation about 
the axis through~$z$ by the angle~$\theta$.  
Consider the diffeomorphism $\psi:M\to M$ defined by
$$
\psi(\theta_1,\theta_2,z_1,z_2) 
:= (\theta_1,\theta_2,z_1,\phi_{z_1,\theta_1}(z_2)),\qquad
\theta_1,\theta_2\in S^1,\qquad z_1,z_2\in S^2.
$$
It acts as the identity on cohomology.
Let $\om\in\Om^2(M)$ be the product symplectic 
form, where both $S^2$ factors have the same area, 
and denote 
$$
a:=[\om]=[\psi^*\om]\in H^2(M;\R).
$$   
McDuff's theorem asserts that $\om$ and $\psi^*\om$ 
can be joined by a path of symplectic 
forms, but not by a path of cohomologous symplectic forms.
The proof that $\om$ and $\psi^*\om$ are not isotopic cannot 
be based on the Gromov--Witten invariants because these are 
invariant under deformation of symplectic forms 
(equivalence relation~(B) in~\ref{para:equiv}).
McDuff's proof does involve the moduli space of holomorphic 
spheres.  The relevant evaluation maps represent the same
homology class but are not homotopic; they have different Hopf 
invariants.  The argument breaks down for symplectic forms
where the $S^2$ factors have different areas, because
in that case the relevant moduli spaces are no longer compact.

\smallskip\noindent{\bf Conclusion.}
{\it For the six-manifold $M=\T^2\times S^2\times S^2$ question~1
has a negative answer and questions~2 and~3 are open problems.}
\end{example}

\begin{example}[{\bf The projective plane}]\label{ex:CP2}\rm
A theorem of Taubes~\cite{TAUBES3,TAUBES4,TAUBES5} asserts
that any two symplectic forms on $M=\CP^2$
with the same volume are diffeomorphic.
The proof uses a theorem of Gromov~\cite{GROMOV1} which requires, 
as an additional hypothesis, the existence of a symplectically
embedded two-sphere (see also~\cite[Theorem~9.4.1]{MS2}).  
The existence of the required symplectic sphere follows from
Taubes' {\it ``Seiberg--Witten equals Gromov''} theorem~\cite{TAUBES5}.   
Combining Taubes' theorem with Moser isotopy one finds that,
for every cohomology class $a\in H^2(M;\R)$ with $a^2\ne0$,
there is a bijection 
$$
\cS_a \cong \Diff(M,a)/\Diff(M,\om),
$$
where $\om$ is the Fubini--Study form re\-pre\-sen\-ting the 
class $a$ and 
\begin{equation*}
\begin{split}
\Diff(M,a) 
&:= \left\{\phi\in\Diff(M)\,|\,\phi^*a=a\right\} \\
&\,= \left\{\phi\in\Diff(M)\,|\,
\phi_*=\id:H_*(\CP^2;\Z)\to H_*(\CP^2;\Z)\right\}.
\end{split}
\end{equation*}
Another theorem of Gromov~\cite{GROMOV1} asserts 
that $\Diff(M,\om)$ retracts onto the isometry group $\PU(3)$ 
of $\CP^2$ (see also~\cite[Theorem~9.5.3]{MS2}).  

\smallskip\noindent
{\bf Conclusion.} {\it Questions~2 and~3 have positive answers. 
Moreover, $\cS_a$ is connected if and only if every diffeomorphism 
of $\CP^2$ that acts as the identity on homology 
is isotopic to the identity}.  
\end{example}

\begin{example}[{\bf The product $S^2\times S^2$}]\label{ex:S2S2}\rm
The discussion of Example~\ref{ex:CP2} carries over 
to $M:= S^2\times S^2$ as follows.  First, every class 
$a\in H^2(M;\R)$ with $a^2\ne 0$ is represented by a symplectic form.  
Second, theorems of Gromov~\cite{GROMOV1} and McDuff~\cite{MCDUFF1}
assert that every symplectic form for which the homology classes 
$A:=[S^2\times\{\mathrm{pt}\}]$ and $B:=[\{\mathrm{pt}\}\times S^2]$ 
(with either orientation) are represented by symplectically embedded 
spheres is diffeomorphic to a standard form (see~\cite[Theorem~9.4.7]{MS2}).  
Third, Taubes' theorem~\cite{TAUBES5} 
establishes the existence of the required symplectic spheres. 
Fourth, a theorem of Gromov~\cite{GROMOV1} asserts that 
the group of symplectomorphisms that preserve $A,B$ 
retracts onto the isometry group $\SO(3)\times\SO(3)$ when
$A,B$ have the same area (see~\cite[Theorem~9.5.1]{MS2}).   
Fifth, a theorem of Abreu and McDuff~\cite{A,AM}
asserts that the symplectomorphism group is connected 
when $A,B$ have different areas.  
Thus $\cS_a\cong\Diff(M,a)/\Diff(M,\om)$ and $\Diff(M,\om)$ has 
two connected components when $\inner{a}{A}=\inner{a}{B}$, 
and is connected otherwise.  

\smallskip\noindent
{\bf Conclusion.}  {\it Questions~2 and~3 have positive answers. 
Moreover, $\cS_a$ is connected if and only if every diffeomorphism 
of ${S^2\times S^2}$ that acts as the identity on homology 
is isotopic to the identity}. 
\end{example}

\begin{example}[{\bf Ruled surfaces}]\label{ex:ruled}\rm
Let $M$ be an orientable smooth four-manifold that admits the 
structure of a fibration over a closed orientable surface $\Sigma$
of positive genus with fibers diffeomorphic to the $2$-sphere:
$$
S^2\longhookrightarrow M \stackrel{\pi}{\longrightarrow} \Sigma.
$$
Fix an orientation of $M$ and an orientation of the fibers.
Let $F\in H_2(M;\Z)$ be the homology class of the fiber.
Call a symplectic form $\om$ on $M$ {\bf compatible with the fibration} 
if it restricts to a symplectic form on each fiber. Call it {\bf compatible 
with the orientations} if its cohomology class $a$ satisfies
$$
a^2>0,\qquad \inner{a}{F}>0.
$$
Here are some basic facts.  

\smallskip\noindent{\bf 1.}
By Seiberg--Witten theory~\cite{TAUBES1,TAUBES2} a cohomology 
class $a\in H^2(M;\R)$ is represented by a symplectic form 
if and only if $a^2\ne0$ and $\inner{a}{F}\ne 0$.  Any such cohomology 
class is uniquely determined by the numbers $a^2$ and $\inner{a}{F}$.

\smallskip\noindent{\bf 2.}
A theorem of McDuff~\cite{MCDUFF1}  (see also~\cite[Theorem~9.4.1]{MS2})
shows that every symplectic form on $M$ that admits a symplectically 
embedded two-sphere in the homology class $F$ or $-F$ is diffeomorphic 
to one that is compatible with the fibration.  
The existence of the required symplectic sphere follows 
from Taubes--Seiberg--Witten theory~\cite{TAUBES5}.

\bigbreak

\smallskip\noindent{\bf 3.}
$M$ admits an orientation preserving diffeomorphism that preserves 
the fibration and reverses the orientation of the fiber, and  
an orientation reversing diffeomorphisms that preserves 
the fibration and the orientation of the fiber.
Thus every symplectic form on $M$ is diffeomorphic to one 
that is compatible with the fibration and orientations.

\smallskip\noindent{\bf 4.}
A theorem of Lalonde--McDuff~\cite{LM1,LM2} asserts that any two symplectic 
forms $\om_0,\om_1$ on $M$ that are compatible with the fibration 
and orientations can be joined by a path of symplectic forms.
They also proved that the path can be chosen in the same 
cohomology class when $[\om_0]=[\om_1]$.  Thus
\begin{equation}\label{eq:LM}
\om_0\stackrel{(\mbox{\tiny{a}})}{\sim}\om_1\qquad\iff\qquad
\om_0\stackrel{(\mbox{\tiny{b}})}{\sim}\om_1\;\;\mbox{ and }\;\;[\om_0]=[\om_1]
\end{equation}
Here~(a) and~(b) denote the equivalence relations in~\ref{para:equiv}.

\smallskip\noindent
{\bf Conclusion.} {\it On a ruled surface any two 
symplectic forms are deformation equivalent, 
and diffeomorphic if they represent the same cohomology class.
In the latter case they are homotopic if and only if they are isotopic.}
\end{example}

\begin{example}[{\bf The one point blow up of the projective plane}]
\label{ex:CP21}\rm
Let
$$
M := \CP^2\#\overline{\CP}^2.
$$
Let $L\in H_2(M;\Z)$ be the homology class of the line in $\CP^2$ 
and $E\in H_2(M;\Z)$ be the class of the exceptional divisor,
both with their complex orientations. They have self-intersection
numbers $L\cdot L=1$ and $E\cdot E=-1$.
Then $M$ admits an orientation reversing diffeomorphism 
that interchanges $L$ and $E$, an orientation preserving 
diffeomorphism that preserves $L$ and reverses $E$, 
and an orientation preserving diffeomorphism that 
reverses $L$ and preserves~$E$.  
By Taubes' theorem~\cite{TAUBES3,TAUBES4,TAUBES5} 
$L$ or $-L$, and $E$ or $-E$, is represented by a symplectically 
embedded sphere for every symplectic form on $M$.  Hence 
a class $a\in H^2(M;\R)$ is represented by a symplectic form 
if and only if  
\begin{equation}\label{eq:LE}
a^2\ne0,\qquad\inner{a}{L}\ne 0,\qquad \inner{a}{E}\ne0.
\end{equation}
By the theorems of McDuff~\cite{MCDUFF1,MCDUFF4} and 
Taubes~\cite{TAUBES5} any two symplectic forms representing 
the same cohomology class are diffeomorphic, as in Example~\ref{ex:ruled}. 
Thus $\cS_a\cong\Diff(M,a)/\Diff(M,\om)$. By~\eqref{eq:LE} 
every diffeomorphism preserving $a$ acts as the identity on homology.  
Moreover, a theorem of Abreu and McDuff~\cite{A,AM} asserts that 
$\Diff(M,\om)$ is connected.  

\smallbreak\noindent
{\bf Conclusion.} {\it Questions~2 and~3 have positive answers.
Moreover, $\cS_a$ is connected if and only if every diffeomorphism of $M$ 
that acts as the identity on homology is smoothly isotopic to the identity.} 
\end{example}

\bigbreak

\begin{example}[{\bf The two point blow up of the projective plane}]
\label{ex:CP22}\rm
Consider the oriented four-manifold
$
M := \CP^2\#\overline{\CP}^2\#\overline{\CP}^2.
$
Let $L\in H_2(M;\Z)$ be the homology class of the line with 
self-intersection number $L\cdot L=1$ and $E_1,E_2\in H_2(M;\Z)$ 
be the homology classes of the exceptional divisors with 
self-intersection numbers $E_i\cdot E_i=-1$. 
Then, by~\eqref{eq:hirzebruch} with $2\chi+3\sigma=7$,
$$
\cC = \left\{c=\PD(nL-n_1E_1-n_2E_2)\,\bigg|\,
\begin{array}{l}
n,\,n_1,\,n_2\mbox{ are odd},\\
n^2-n_1^2-n_2^2=7
\end{array}\right\}.
$$
The set $\cE$ of all classes $E\in H_2(M;\Z)$
satisfying $E\cdot E=-1$ is given by
$$
\cE = \{E=mL+m_1E_1+m_2E_2\,|\,m^2-m_1^2-m_2^2=-1\}.
$$
A theorem of Li--Liu~\cite{LL1,LL2} asserts the 
following for every $a\in H^2(M;\R)$.  

\smallskip\noindent{\bf (I)}
{\it $\cC_\symp=\cC$.}

\smallskip\noindent{\bf (II)}
{\it The diffeomorphism group acts transitively on $\cC_\symp$.}

\smallskip\noindent{\bf (III)}
{\it $\cS_a\ne\emptyset$ if and only if $a^2>0$ 
and $\inner{a}{E}\ne 0$ for every $E\in\cE$.}

\smallskip\noindent{\bf (IV)}
{\it $\cS_a\ne\emptyset\implies\#\cC_a=1$.}

\smallskip\noindent{\bf (V)}
{\it If $\om_0,\om_1$ are symplectic forms on~$M$,
then there exists a diffeomorphism~$\phi$ of~$M$ 
such that $\om_0$ and $\phi^*\om_1$ can be joined 
by a path of symplectic forms.}

\smallskip\noindent
Assuming Seiberg--Witten theory, assertion~(VI) below is Lemma~3.11
in Karshon--Kessler--Pinsonnault~\cite{KKP} and~(VII) 
is Theorem~1.1 in McDuff~\cite{MCDUFF5}.

\smallskip\noindent{\bf (VI)}
{\it Two symplectic forms $\om_0,\om_1$ on~$M$ have the 
same first Chern class if and only if there exists a 
diffeomorphism~$\phi$ of~$M$, inducing the identity on homo\-logy, 
such that $\om_0$ and $\phi^*\om_1$ can be joined 
by a path of symplectic forms.}

\smallskip\noindent{\bf (VII)}
{\it Two symplectic forms $\om_0,\om_1$ on~$M$ represent the
same cohomology class if and only if there exists 
a diffeomorphism $\phi$ of $M$, inducing the identity 
on homology, such that $\phi^*\om_1=\om_0$.}

\smallskip
In particular, questions~2 and~3 have positive answers.
The proofs rely on the following observations. 

\smallskip\noindent{\bf 1.}
For every $c\in\cC$ there exist precisely three 
classes $E^c_1,E^c_2,E^c_3\in\cE$ such that 
$\inner{c}{E^c_i} = 1$.  The numbering can be chosen such that
\begin{equation}\label{eq:Ec}
E^c_1\cdot E^c_3 = E^c_2\cdot E^c_3=1,\qquad 
E^c_1\cdot E^c_2 = 0.
\end{equation}
Each class $E^c_i$ is represented by an embedded sphere.
For the standard first Chern class 
$c_\standard=\PD(3L-E_1-E_2)\in\cC$
these can be chosen as the classes 
$E^{c_\standard}_1=E_1$, $E^{c_\standard}_2=E_2$, 
$E^{c_\standard}_3=L-E_1-E_2$.
For general elements $c\in\cC$ this follows from the fact that,
by the results of Wall~\cite{WALL1,WALL2,WALL3,WALL4}, the 
diffeomorphism group acts transitively on $\cC$. 
(See Li--Liu~\cite{LL1,LL2} for details.) 

\bigbreak

\smallskip\noindent{\bf 2.}
By Taubes--Seiberg--Witten theory~\cite{TAUBES3,TAUBES4,TAUBES5}
the homology classes $E^c_1,E^c_2,E^c_3$ are represented by 
symplectically embedded spheres for every symplectic form 
on $M$ with first Chern class $c$.  

\smallskip\noindent{\bf 3.}
A cohomology class $a\in H^2(M;\R)$ can be represented by a 
symplectic form $\om$ with first Chern class $c_1(\om)=c$ 
if and only if for every $E\in H_2(M;\Z)$
\begin{equation}\label{eq:a-symp}
E\cdot E=-1,\quad \inner{c}{E}=1
\qquad\implies\qquad \inner{a}{E}>0.
\end{equation}
By~2.\ condition~\eqref{eq:a-symp} is necessary for the 
existence of~$\om$.  The converse follows from the 
Nakai--Moishezon criterion. (See Demazure~\cite{DEMAZURE}
and also McDuff--Polterovich~\cite{MP} 
for a direct symplectic proof.) 

\smallskip\noindent{\bf 4.}
Let $a\in H^2(M;\R)$ such that $a^2>0$. Close examination 
shows that there exists a cohomology class $c\in\cC$
such that condition~\eqref{eq:a-symp} holds
if and only if $\inner{a}{E}\ne0$ for {\it all} $E\in\cE$. 
(See Li--Liu~\cite{LL1,LL2} for details.) 

\smallskip\noindent
Assertions~(I), (II), (III) follow from 1.-4.\ 
and~(IV) follows from Corollary~A.

\begin{figure}[htp] 
\centering 
\includegraphics[scale=0.35]{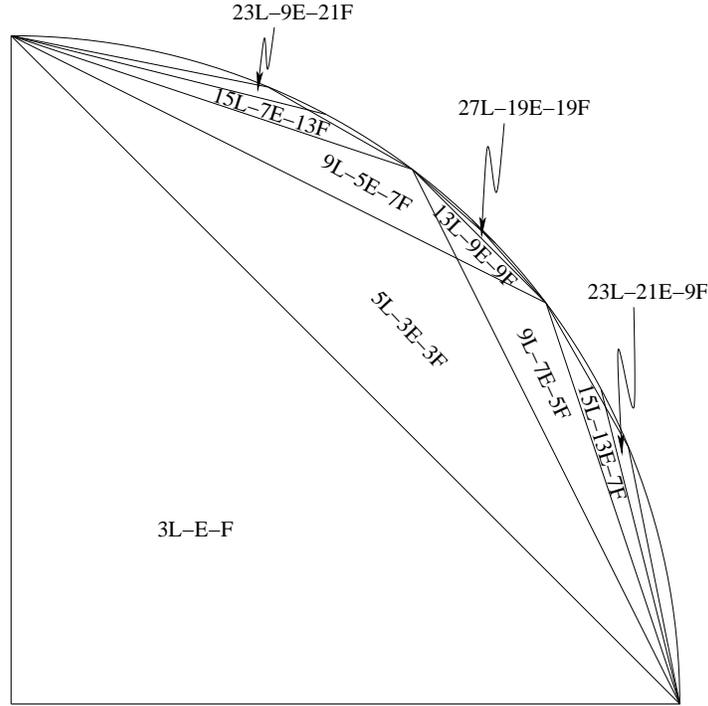}      
\caption{{The symplectic chamber structure on 
$\CP^2\#2\overline{\CP}^2$.}}\label{fig:CP22}      
\end{figure} 

\smallskip\noindent{\bf 5.}
To understand~(IV) geometrically, it is convenient 
to slightly change the point of view.
Fix a cohomology class $c\in\cC$ and define 
\begin{equation}\label{eq:sympcone}
\begin{split}
\cE^c 
&:= \left\{E\in H_2(M;\Z)\,|\,
E\cdot E=-1,\,\inner{c}{E}=1\right\},\\
\cK^c
&:= 
\left\{a\in H^2(M;\R)\,|\,a^2>0,\,
\inner{a}{E}>0\;\forall\,E\in\cE^c\right\},\\
\cS^c 
&:= 
\left\{[\rho]\in H^2(M;\R)\,|\,\rho\in\Om^2(M),\,
d\rho=0,\,\rho\wedge\rho>0,\,c_1(\rho)=c\right\}.
\end{split}
\end{equation}
Then $\cS^c=\cK^c$ (see 3.\ above).
These cones define a chamber structure on $H^2(M;\R)$.   
(See Figure~\ref{fig:CP22}, where the standard basis 
of $H_2(M;\Z)$ is denoted by $L,E,F$ and the symplectic 
cohomology classes satisfy the normalization conditions 
$\inner{a}{L}=1$, $\inner{a}{E}>0$, $\inner{a}{F}>0$.)
The chamber structure is determined by straight lines
connecting pairs of rational points on the unit circle.
In this notation~(IV) asserts that the chambers are 
disjoint.

\smallskip\noindent{\bf 6.}
Here is a sketch of a proof of~(V) and~(VI),
which was worked out in a discussion with Yael Karshon.
Construct a model for the symplectic blowup 
$M=\CP^2\#2\overline{\CP}^2$ as follows.  
Choose distinct points $p_1,p_2\in\CP^2\setminus\CP^1$
and disjoint neighborhoods $U_1,U_2\subset\CP^2\setminus\CP^1$
of $p_1,p_2$, respectively, equipped with holomorphic 
coordinate charts 
$$
\phi_i:B^4_\eps\to U_i,
$$
defined on the $\eps$ ball in $\R^4=\C^2$, centered at the origin.
Then use the standard complex blowup construction 
at $p_1$ and $p_2$ to define $M$ as a complex  manifold.  
There is a holomorphic projection 
$$
\pi:M\to\CP^2
$$ 
with singular values $p_1,p_2$.  Their preimages 
$S_1:=\pi^{-1}(p_1)$ and $S_2:=\pi^{-1}(p_2)$ are 
the exceptional divisors that, with their complex orientations,
represent the classes $E_1,E_2$.  
The projection $\pi$ is a holomorphic diffeomorphism from 
$M\setminus(S_1\cup S_2)$ to $\CP^2\setminus\{p_1,p_2\}$
and the holomorphic sphere $S:=\pi^{-1}(\CP^1)$ represents 
the class $L$.  

To define a symplectic form on $M$, let $\om_\FS$ denote the 
Fubini--Study form on $\CP^2$ such that $\CP^1$ has area $\pi$.
Choose an exact perturbation $\om_0$ of $\om_\FS$ such that 
$\phi_i^*\om_0$ is standard on $B^4_\eps$ for $i=1,2$.  
(Warning: $\om_0$ cannot be equal to $\om_\FS$.)
Now use the construction in~\cite[Lemma~7.15]{MS1} 
to obtain a symplectic form $\rho_0$ on~$M$ such that 
$\rho_0=\pi^*\om_0$ on $M\setminus\pi^{-1}(U_1\cup U_2)$ 
and the exceptional spheres $S_1,S_2$ are symplectic 
submanifolds of small areas.

Let $\om$ be any symplectic form on~$M$ with first Chern class
$c:=c_1(\om)$ and choose $E^c_1,E^c_2,E^c_3$ as in part~1 above
such that~\eqref{eq:Ec} holds.  Denote
$$
L^c := E^c_1+E^c_2+E^c_3,\qquad
\lambda_i := \inner{[\om]}{E^c_i},\qquad i=1,2,3,
$$
and assume $\inner{[\om]}{L^c}=\pi$.
By Taubes--Seiberg--Witten theory~\cite{TAUBES5}, 
the classes $E^c_1,E^c_2,L^c$ are represented 
by pairwise disjoint symplectically embedded spheres 
$C_1,C_2,C$. Hence there exist disjoint
Weinstein neighborhoods $N_1,N_2\subset M\setminus C$ 
of $C_1,C_2$ in which the symplectic form 
is standard (see~\cite[Theorem~3.30]{MS1}). 
The symplectic blowdown construction replaces
$C_1,C_2$ by closed four-balls of radii $r_1,r_2$ 
with $\pi r_i^2=\lambda_i$.  This results in a symplectic 
manifold $(\widehat{M},\widehat{\om})$,
equipped with symplectic embeddings 
$\widehat{\psi}_i:B^4_{r_i+\delta}\to\widehat{M}$ 
of closed four-balls with images
$\widehat{N}_i:=\widehat{\psi}_i(B^4_{r_i+\delta})$,
and a symplectomorphism
$$
\widehat{\pi}:(M\setminus(N_1\cup N_2),\om)\to 
(\widehat{M}\setminus(\widehat{N}_1\cup\widehat{N}_2),\widehat{\om}).
$$
Thus $C$ descends to a symplectically embedded sphere
$\widehat{C}:=\widehat{\pi}(C)\subset\widehat{M}$, disjoint from 
$\widehat{N}_1\cup\widehat{N}_2$, with area $\pi$ and 
self-intersection number $\widehat{C}\cdot\widehat{C}=1$. 
The homology class of~$\widehat{C}$ generates $H_2(\widehat{M};\Z)$. 
Hence a theorem of Gromov~\cite{GROMOV1}
and McDuff~\cite{MCDUFF1} (see also~\cite[Theorem~9.4.1]{MS2}) 
asserts that there exists a symplectomorphism
$g:(\widehat{M},\widehat{C},\widehat{\om})\to(\CP^2,\CP^1,\om_0)$. 
Consider the embeddings 
$
\psi_i:=g\circ\widehat{\psi}_i:B^4_{r_i+\delta}\to\CP^2.
$ 
Composing $g$ with a suitable Hamiltonian isotopy
we may assume without loss of generality that 
$\psi_i(0)=p_i$ and that $\psi_i$ agrees with 
$\phi_i$ near the origin. Thus 
$g(\widehat{N}_i) = \psi_i(B^4_{r_i+\delta})$
is a neighborhood of $p_i$ and the complement of their
union is symplectomorphic to a subset of $M$ under~$\pi^{-1}$.  
The formula in the proof of~\cite[Lemma~7.15]{MS1} gives 
rise to a diffeomorphism $f:M\to M$ such that
$f(C_i) = S_i = \pi^{-1}(p_i)$, $f(N_i) = \pi^{-1}(g(\widehat{N}_i))$
and $f|_{M\setminus(N_1\cup N_2)}=\pi^{-1}\circ g\circ\widehat{\pi}$.
The diffeomorphism $f$ satisfies
$$
f_*L^c = L,\qquad f_*E^c_1 = E_1,\qquad f_*E^c_2 = E_2.
$$
The symplectic form $\rho:=(f^{-1})^*\om$ agrees with $\rho_0$ 
on $M\setminus\pi^{-1}(g(\widehat{N}_1)\cup g(\widehat{N}_2))$.
It can be joined to $\rho_0$ by a path of symplectic forms via a 
scaling argument. If $c_1(\om)=c_\standard$ one can choose 
$E^c_i=E_i$ and then $f$ induces the identity on homology.
This proves~(V) and~(VI).  

\smallskip\noindent{\bf 7.}
The proof of~(VII) hinges on the fact that~\eqref{eq:LM} holds for the 
two point blowup of $\CP^2$ and symplectic forms
with the standard first Chern class (i.e.\ two cohomologous
symplectic forms with first Chern class $c_\standard$
are homotopic if and only if they are isotopic).  
This was proved by McDuff in~\cite{MCDUFF4,MCDUFF5}.
Thus~(VII) follows from~(VI) and Corollary~A.
\end{example}

\begin{example}[{\bf Del Pezzo surfaces}]
\label{ex:delPezzo}\rm
The discussion in Example~\ref{ex:CP22} carries over
almost verbatim to all del Pezzo surfaces
$$
M := \CP^2\#k\overline{\CP}^2,\qquad
1\le k\le 8.
$$

\smallskip\noindent{\bf 1.}
The characterization of the symplectic cone $\cS^c=\cK^c$
in~\eqref{eq:sympcone} remains unchanged and the number 
of exceptional classes is
$$
\#\left\{E\in H_2(M;\Z)\,\bigg|\,\begin{array}{l}
E\cdot E=-1,\\
\inner{c}{E}=1
\end{array}
\right\}
= \left\{\begin{array}{rl}
1,&\mbox{if }k=1,\\
3,&\mbox{if }k=2,\\
6,&\mbox{if }k=3,\\
10,&\mbox{if }k=4,\\
16,&\mbox{if }k=5,\\
27,&\mbox{if }k=6,\\
56,&\mbox{if }k=7,\\
240,&\mbox{if }k=8.
\end{array}\right.
$$
(For $k=6$ these are the 27 lines on a cubic in $\CP^3$.)
Thus condition~\eqref{eq:a-symp} is still necessary
and sufficient for the existence of a symplectic
form that represents the class $a$ 
and has first Chern class $c$. That the condition 
is necessary follows again from Taubes--Seiberg--Witten 
theory~\cite{TAUBES5}.  That the condition is also sufficient
for the standard class 
$$
c_\standard:=\PD\bigl(3L-E_1-\cdots-E_k\bigr)
$$
is a classical result in K\"ahler geometry.
That the condition is sufficient in general, was proved 
by Li--Liu~\cite{LL1,LL2}, by reducing the assertion 
to the standard chamber via the results of 
Wall~\cite{WALL1,WALL2,WALL3}. 

\smallskip\noindent{\bf 2.}
Assertion~\eqref{eq:LM} continues to hold for all blowups 
of $\CP^2$, i.e.\ two cohomologous symplectic forms 
are homotopic if and only if they are isotopic.
This was first proved by McDuff~\cite{MCDUFF4,MCDUFF5} for $k\le2$,
then by Biran~\cite{BIRAN1} for $3\le k\le 6$,
then by McDuff~\cite{MCDUFF2} for all~$k$.
The proofs (for $k\ge 3$) are based on the inflation 
techniques of Lalonde--McDuff~\cite{LALONDE,LM1,LM2,MCDUFF2}. 

\smallskip\noindent{\bf 3.}
With this understood it follows by the same arguments as 
in Example~\ref{ex:CP22}, parts~6 and~7,
that assertions~(V), (VI), (VII) continue to hold for all 
del Pezzo surfaces. (See Li--Liu~\cite{LL1,LL2} 
and Karshon--Kessler--Pinsonnault~\cite{KKP}.)

\smallskip\noindent{\bf Conclusion.}
{\it Assertions~(I)-(VII) in~\ref{ex:CP22} remain valid 
for all del Pezzo surfaces.  Thus questions~2 
and~3 have positive answers and $\cC_\symp=\cC$. 
Two cohomologous symplectic forms are homotopic 
if and only if they are isotopic.} 
\end{example}

\begin{example}[{\bf Higher blowups of the projective plane}]
\label{ex:elliptic}\rm
Consider the closed oriented smooth four-manifolds
$M := \CP^2\#k\overline{\CP}^2$, $k\ge 9$.

\smallskip\noindent{\bf 1.}~When $k\ge 9$ 
it is still true that a class $a\in H^2(M;\R)$ 
is represented by a symplectic form with first Chern 
class $c=c_\standard$ if and only if ${a^2>0}$ and~\eqref{eq:a-symp} 
holds. Thus, in the notation of~\eqref{eq:sympcone}, 
the symplectic cone is given by
$\cS^{c_\standard} = \cK^{c_\standard}$.
For $k>9$ this no longer follows from K\"ahler geometry.
It is an existence theorem in symplectic topology by
Li--Liu~\cite{LL2}. An earlier theorem by 
Biran~\cite{BIRAN2} asserts that 
$\cK^{c_\standard}\subset\cS^{c_\standard}
\subset\overline{\cK}^{c_\standard}$.
The proofs are based on Taubes--Seiberg--Witten 
theory and the inflation techniques developed 
by Lalonde and McDuff~\cite{LALONDE,LM1,LM2,MCDUFF5,MCDUFF2}.
The Nagata conjecture asserts that the 
K\"ahler cone for the standard complex structure 
should agree with the symplectic cone.
(See Biran~\cite{BIRAN4} for a 
symplectic approach to this question.)

\smallskip\noindent{\bf 2.}~Consider 
the {\bf elliptic surface}
$
M := E(1) := \CP^2\#9\overline{\CP}^2.
$
This manifold admits the structure of a holomorphic 
Lefschetz fibration over $\CP^1$ with elliptic curves
as regular fibers and twelve singular fibers. 
(Choose two cubics in general position and
blow up their nine points of intersection.)
In the standard basis $L,E_1,\dots,E_9$ of $H_2(M;\Z)$
the homology class of the fiber is
$F := 3L-E_1-E_2-\cdots-E_9$. Its Poincar\'e dual 
$c_\standard = \PD(F) = c_1(TM,J_\standard)$
is the first Chern class of the standard complex structure.

\smallskip\noindent{\bf 3.}~For $k\ge 9$ the 
diffeomorphism group no longer 
acts transitively on~$\cC$, however, it acts 
transitively on $\cC_\symp\subsetneq\cC$.  
For $k=9$ the theory of Wall~\cite{WALL3} 
is still applicable and shows that, for every $c\in\cC$, 
there is a unique odd integer $m\ge1$ such that $c$ is diffeomorphic 
to $m c_\standard$. By Taubes' results, $c$ is the first 
Chern class of a symplectic form if and only if ${m=1}$.  
For $k\ge9$, Li--Liu proved 
in~\cite[Theorem~1]{LL2} that 
\begin{equation}\label{eq:Csymp}
\cC_\symp = \{c\in\cC\,|\,
\exists \phi\in\Diff(M)\mbox{ s.t. }
\phi^*c=c_\standard\}.
\end{equation}
In~\cite[Theorem~D]{LL1} they showed
that any two symplectic forms with the standard
first Chern class are deformation equivalent.
Assuming~\eqref{eq:Csymp}, the proof of~(V)-(VII) for 
all $k$ follows the same line of argument as in 
Example~\ref{ex:CP22}, parts~6 and~7, in the case $k=2$.
(See Karshon--Kessler--Pinsonnault~\cite{KKP}.)

\smallskip\noindent{\bf 4.}
There is a chamber structure on $H^2(M;\R)$ 
as in Examples~\ref{ex:CP22} and~\ref{ex:delPezzo}.
However, the set $\cE^c$ is now infinite 
for each $c\in\cC_\symp$.  For more details see 
Biran~\cite{BIRAN2,BIRAN3,BIRAN4} and Li-Liu~\cite{LL1,LL2}.

\smallskip\noindent{\bf Conclusion.}
{\it Assertions~(II)-(VII) in~\ref{ex:CP22}
remain valid for $k\ge 9$. Thus questions~2 and~3 
have positive answers. However, $\cC_\symp\subsetneq\cC$.
Two cohomologous symplectic forms are homotopic
if and only if they are isotopic.}
\end{example}

\begin{example}[{\bf Deformation equivalence}]
\label{ex:D}\rm\hfill

\smallskip\noindent{\bf 1.}
This example is due to Smith~\cite{SMITH1}.
Consider the four-manifolds
$$
X := E(5) := E(1) \#_{\T^2} E(1) 
\#_{\T^2} E(1) \#_{\T^2} E(1) \#_{\T^2} E(1),\quad
Y := Z_5 \# 5\overline{\CP}^2,
$$
where $Z_5\subset\CP^3$ is a degree five hypersurface.
Then $X$ admits the structure of a Lefschetz fibration
with generic fiber $F$ a torus and 60 singular fibers. 
Both manifolds are simply connected, 
not spin, and have the Betti numbers
$$
b^+=9,\qquad b^-=49,\qquad \sigma=-40,\qquad \chi=60.
$$
Hence they are homeomorphic and hence
$X\times X$ and $Y\times Y$ are diffeomorphic.  
However, the first Chern class of $X$ with its 
standard K\"ahler structure is divisible by three (it is $-3\PD(F)$), 
while the first Chern class of $Y$ is primi\-tive. These 
divisibility properties carry over to the product.   
Hence there exist symplectic forms 
${\om_X\in\Om^2(X\times X)}$ and
${\om_Y\in\Om^2(Y\times Y)}$ on diffeomorphic
manifolds such that 
$$
c_1(\om_X)\ne\phi^*c_1(\om_Y)
$$
for every diffeomorphism $\phi:X\times X\to Y\times Y$.
This gives rise to two symplectic forms on the same
manifold that are not even related by the weakest 
equivalence relation~(D) in~\ref{para:equiv}.
Hence they are not deformation equivalent.

\smallskip\noindent{\bf 2.}
In~\cite{SMITH1,SMITH2} Smith constructed a simply 
connected four-manifold $X$ that admits two symplectic 
forms $\om_0$, $\om_1$ such that $c_1(\om_0)$ is divisible 
by three and $c_1(\om_1)$ is primitive.  His four-manifold
is obtained by forming a fiber connected sum of $\T^4$ 
with five copies of $E(1)$. As in 1.\ above the first Chern classes
of $\om_0$ and $\om_1$ are not diffeomorphic, and hence
$\om_0$ and $\om_1$ are not deformation equivalent.
In fact, for each integer $N$ he constructed a simply 
connected four-manifold $X_N$ with at least $N$ different
deformation equivalence classes of symplectic forms,
distinguished by the divisibility properties of their first
Chern classes.  Taking the product $X_N\times\T^{2n}$, 
one obtains a $(4+2n)$-manifold with $N$ pairwise 
deformation inequivalent symplectic forms.

\smallskip\noindent{\bf 3.}
An earlier example of a symplectic four-manifold 
with two deformation inequivalent symplectic forms
was constructed by McMullen--Taubes~\cite{MT}.  
In their example the symplectic forms are distinguished by their 
Seiberg--Witten invariants (see also~\cite{VIDUSSI2}). 
In~\cite{VIDUSSI1} Vidussi constructed homotopy K3's 
with arbitrarily many deformation inequivalent 
symplectic forms.

\smallskip\noindent{\bf 4.}
In all these examples it is not clear whether  
two deformation inequivalent symplectic forms 
can be found in the same cohomology class. 
\end{example}

\begin{example}[{\bf Dehn--Seidel twists}]\label{ex:seidel}\rm
In his PhD thesis Seidel constructed symplectomorphisms 
on many symplectic four-manifolds $(M,\om)$ that 
are smoothly, but not symplectically, isotopic 
to the identity.  He proved the following two theorems
(see~\cite{SEIDEL1} and~\cite[Cor~2.9 \& Thm~0.5]{SEIDEL3}).

\medskip\noindent{\bf Theorem~A.}
{\it Let $M$ be a closed symplectic four-manifold with Betti numbers
$$
b_1:=\dim\,H^1(M;\R)=0,\qquad
b_2:=\dim\,H^2(M;\R)\ge 3.
$$
Assume $M$ is minimal (i.e.\ it does not contain a symplectically 
embedded two-sphere with self-intersection number minus one).

If $\Lambda\subset M$ is a Lagrangian sphere, then the square 
$$
\phi =\tau_\Lambda\circ\tau_\Lambda
$$
of the Dehn--Seidel twist determined by $\Lambda$
is smoothly, but not symplectically, isotopic to the identity.}

\medskip\noindent{\bf Theorem~B.}
{\it Let $M$ be an algebraic surface and a complete intersection, 
not diffeomorphic to $\CP^2$ or $\CP^1\times\CP^1$.  
Then there is a symplectomorphism $\phi:M\to M$ that is smoothly, 
but not symplectically, isotopic to the identity.}

\medskip\noindent
Here are some remarks.  For a much more detailed and wide ranging
discussion of examples and ramifications of these theorems 
see Seidel~\cite{SEIDEL3}. 

\smallskip\noindent{\bf 1.}
The assumptions of Theorem~A allow for examples 
of symplectic four-manifolds that do not admit K\"ahler structures
(see for example Gompf--Mrowka~\cite{GM} 
and Fintushel--Stern~\cite{FS}).
The assumptions of Theorem~B allow for examples that 
are not minimal, such as the six-fold blowup of 
the projective plane (cubics in $\CP^3$).  

\smallskip\noindent{\bf 2.} 
The symplectomorphism $\phi=\tau_\Lambda^2$ in Theorem~A
is smoothly isotopic to the identity by an isotopy localized near~$\Lambda$.  
Seidel computed the Floer cohomology group $\HF^*(\tau_\Lambda)$ 
with its module structure over the quantum cohomology ring,
via his exact sequence~\cite{SEIDEL1,SEIDEL2}.
As a result he was able to show that, under the
assumptions of minimality and $b_2\ge3$, the Floer cohomology 
group $\HF^*(\tau_\Lambda)$ is not isomorphic to the Floer 
homology group 
$$
\HF_*(\tau_\Lambda)=\HF^*(\tau_\Lambda^{-1}).
$$
Hence $\tau_\Lambda$ and $\tau_\Lambda^{-1}$ are not 
Hamiltonian isotopic. When $b_1=0$ it then follows that they
are not symplectically isotopic.  

\smallbreak

\smallskip\noindent{\bf 3.}
The assumption $b_2\ge3$ cannot be removed in Theorem~A.  
An example is $M=\CP^1\times\CP^1$
with its monotone symplectic structures
and with $\Lambda$ equal to the anti-dia\-go\-nal. 
In this example $\phi=\tau_\Lambda^2$ {\it is} symplectically 
isotopic to the identity, by Gromov's theorem in~\cite{GROMOV1} 
(see Example~\ref{ex:S2S2}).

\smallskip\noindent{\bf 4.}
The assumption of minimality cannot be removed in Theorem~A.
Examples are the blowups $\CP^2\#k\overline{\CP}^2$ with $2\le k\le4$.
On these manifolds there exist symplectic forms 
(monotone in the cases $k=3,4$) and Lagrangian spheres
such that the squares of the Dehn--Seidel 
twists {\it are} symplectically isotopic to the identity 
(see~\cite[Examples~1.10 and~1.12]{SEIDEL3}).  
In contrast, for $5\le k\le8$, the square of a Dehn--Seidel twist 
in the $k$-fold blowup of the projective plane with its 
monotone symplectic form is {\it never} symplectically
isotopic to the identity (see~\cite[Example~2.10]{SEIDEL3}).

\smallskip\noindent{\bf 5.}
Seidel showed that, by an arbitrarily small perturbation of the 
cohomology class of $\om$, the square of the Dehn--Seidel twist
deforms to a symplectomorphism that {\it is} symplectically 
isotopic to the identity.

\smallskip\noindent{\bf 6.}
If $(M,\om)$ is a symplectic four-manifold satisfying 
the assumptions of Theorem~A or Theorem~B then the 
space $\cS_\om$ of all symplectic forms $\rho$ on $M$
that are isotopic to $\om$ is not simply 
connected (see~\ref{para:diffeo}).
However, the nontrivial loops in $\cS_\om$
that arise from Dehn--Seidel twists {\it are} contractible 
in the space of nondegenerate $2$-forms.

\smallskip\noindent{\bf 7.}
In~\cite{KRONHEIMER} Kronheimer used Seiberg--Witten theory 
to prove the existence of symplectic four-manifolds 
$(M,\om)$ such that $\cS_\om$ is not simply connected.
In fact, he developed a method for constructing, 
for each integer $n\ge 1$, symplectic four-manifolds 
$(M,\om)$ such that $\pi_{2n-1}(\cS_\om)\ne 0$.
\end{example}

\begin{example}[{\bf The K3-surface}]\label{ex:K3}\rm
As a smooth manifold the K3-surface can be constructed as
the fiber connected sum 
$$
K3:=E(2)=E(1)\#_{\T^2}E(1).
$$
On the K3-surface with its standard orientation, 
every cohomology class $a\in H^2(K3;\R)$ with $a^2>0$
is represented by a symplectic form (with an associated 
hyperk\"ahler structure).  
By Taubes' results~\cite{TAUBES1,TAUBES2}, 
every symplectic form on the K3-surface 
with the standard orientation 
has first Chern class zero (see~\ref{para:SW}).
Hence it follows from Corollary~B in 
Section~\ref{sec:SYMP} that 
{\it any two symplectic forms on the K3-surface,
compatible with the standard orientation,
are homotopic as nondegenerate $2$-forms.}
\end{example}

\begin{example}[{\bf The four-torus}]\label{ex:torus}\rm
Every cohomology class $a\in H^2(\T^4;\R)$ with 
$a^2\ne 0$ is represented by a symplectic form and,
by Taubes' result in~\cite{TAUBES1}, every symplectic form 
on $\T^4$ has first Chern class zero (see~\ref{para:SW}).  
There are infinitely many homotopy classes 
of nondegenerate $2$-forms with first Chern class zero 
and compatible with a fixed orientation.
The set of such homotopy classes 
is in bijective correspondence to 
the set $\Z/2\Z\times H^3(\T^4;\Z)$ 
(see~\ref{para:homotopyJ}). 
The Conolly-L\'e-Ono argument in~\cite{CLO}
removes only half of these classes as candidates 
for containing a symplectic form 
(see Corollary~B in Section~\ref{sec:SYMP}). 
{\it It remains an open question whether any two
symplectic forms on $\T^4$, that induce the same 
orientation, are homotopic as nondegenerate $2$-forms.}
\end{example}

\begin{example}[{\bf The one point blowup of the four-torus}]
\label{ex:torus1}\rm
A recent theorem by Latschev--McDuff--Schlenk~\cite{LMS} 
asserts that a closed four-ball admits a symplectic 
embedding into the four-torus with any constant coefficient 
symplectic form if and only if the volume of the four-ball
is smaller than the volume of the four-torus. 
This settles the existence problem for the one point 
blowup of the four-torus
$$
M:=\T^4\#\overline{\CP}^2.
$$
Namely, let $E\in H_2(M;\Z)$ be the homology class 
of the exceptional divisor.  Then a cohomology class 
$a\in H^2(M;\R)$ can be represented by a 
symplectic form if and only if 
$$
a^2\ne 0,\qquad \inner{a}{E}\ne 0.
$$
Thus the symplectic cone of $M$ is strictly bigger 
than the K\"ahler cone. The uniqueness problem for 
symplectic forms on $M$ is still far from understood.
\end{example}

\begin{example}[{\bf K\"ahler surfaces of general type}]
\label{ex:generaltype}\rm
Let $(M,\om,J)$ be a minimal K\"ahler surface of general type
with first Chern class $c:=c_1(\om)$ such that 
$$
b^+\ge 2,\qquad c^2 >0,\qquad c\cdot[\om] < 0.
$$
Then $\pm c$ are the only Seiberg--Witten basic classes.
Hence, by Taubes' results in~\cite{TAUBES1},
every symplectic form~$\rho$ on~$M$ satisfies 
$$
c_1(\rho)\cdot[\rho]<0,\qquad c_1(\rho)=\pm c.
$$   
Thus any two symplectic forms, compatible with the orientation,
have the same first Chern class up to sign.
{\it It is an open question whether every cohomology class
$a$ with $a^2>0$ and $a\cdot c<0$ can be represented
by a symplectic form with first Chern class~$c$} 
(as conjectured by Tianjun Li).
\end{example}


\section{Discussion}\label{sec:D}

Here are some of the conclusions that can be drawn from the examples 
in Section~\ref{sec:EX} about the equivalence relations in~\ref{para:equiv}.

\smallskip\noindent{\bf 1.}~Many 
symplectic manifolds (with nonzero first Chern classes) 
admit anti-symplectic involutions.  Thus two symplectic 
forms related by~(A) in~\ref{para:equiv} (a diffeomorphism) 
need not be related by~(d) (same first Chern class).

\smallskip\noindent{\bf 2.}~Example~\ref{ex:McDuff} 
shows that two cohomologous and diffeomorphic 
symplectic forms on a closed manifold $M$
that can be joined by a path of symplectic forms 
need not be isotopic.  Thus
$$
[\om_0]=[\om_1],\quad 
\om_0\stackrel{(\mbox{\tiny{A}})}{\sim}\om_1,\quad
\om_0\stackrel{(\mbox{\tiny{b}})}{\sim}\om_1
\qquad\notimplies\qquad
\om_0\stackrel{(\mbox{\tiny{a}})}{\sim}\om_1.
$$
However, for ruled surfaces and all blowups of the projective plane
it is known that $[\om_0]=[\om_1]$ and~(b) imply~(a)
(see Examples~\ref{ex:CP2}-\ref{ex:elliptic}).

\smallskip\noindent{\bf 3.}~{\it It is an open question 
whether two symplectic forms on a closed manifold 
that can be joined by a path of nondegenerate $2$-forms 
can always be joined by a path of symplectic forms.}  
In other words, it is an open question whether the 
equi\-valence relations~(b) and~(c) agree. 
For open manifolds they agree.

\smallskip\noindent{\bf 4.}~For all closed 
simply connected smooth four-manifolds with 
nonzero Euler characteristic, the equivalence relation~(c)
(homotopic as nondegenerate $2$-forms) 
agrees with~(d) (same first Chern class), 
by Corollary~B in Section~\ref{sec:SYMP}.

\smallskip\noindent{\bf 5.}~For all
blowups of $\CP^2$ the relations~(c)
and~(d) agree with a stronger form of~(B) 
(a diffeomorphism acting as the identity on homology, 
followed by a path of symplectic forms); 
see~\ref{ex:CP22}-\ref{ex:elliptic}.
The symplectic forms on these manifolds 
have only one equivalence class with respect to~(B), (C), (D).

\smallskip\noindent{\bf 6.}~{\it In dimensions $2n\ge6$ it 
is an open question whether two cohomologous symplectic forms 
on a closed manifold always have the same first Chern class.}
For closed oriented smooth four-manifolds
it follows from Seiberg--Witten theory that
$[\om_0]=[\om_1]$ implies $c_1(\om_0)=c_1(\om_1)$
(Corollary~A in Section~\ref{sec:SYMP}). 

\smallskip\noindent{\bf 7.}~{\it 
It is an entirely different question whether 
any two symplectic forms, compatible with the orientation, 
have diffeomorphic first Chern classes~up~to~sign.}
Examples with positive answers include all symplectic 
four-mani\-folds with $b^+=1$ (see Li--Liu~\cite{LL2})
and Examples~\ref{ex:K3}-\ref{ex:generaltype}. 
For negative answers in dimension four and higher 
see Example~\ref{ex:D}.

\smallskip\noindent{\bf 8.}~Example~\ref{ex:D}, part 2,
shows that in any dimension $2n\ge 4$ there exist 
closed manifolds with pairs of symplectic
forms whose first Chern classes are not diffeomorphic
(even up to sign); hence they are not related by~(D).

\bigbreak

\smallskip\noindent{\bf 9.}~{\it It is an open question 
whether there is any closed four-manifold $M$ and 
any cohomology class $a\in H^2(M;\R)$ such that
$\cS_a$ is nonempty and connected.}

\smallskip\noindent{\bf 10.}~{\it It is an open question 
whether there is any closed four-manifold~$M$ and 
any cohomology class $a\in H^2(M;\R)$ such that
$\cS_a$ is disconnected.}

\subsection*{Conjectures}

\begin{conjecture}[{\bf Donaldson's four-six conjecture}]
\label{con:don46}
Let $\sigma$ be a symplectic form on the $2$-sphere.
Let $(X,\om_X)$ and $(Y,\om_Y)$ 
be closed symplectic four-manifolds.  
Then $X$ and $Y$ are diffeomorphic
if and only if the symplectic six-manifolds 
$(X\times S^2,\om_X\times\sigma)$ 
and $(Y\times S^2,\om_Y\times\sigma)$
are deformation equivalent, i.e.\ there exists 
a diffeomorphism $\phi:X\times S^2\to Y\times S^2$ 
such that $\om_X\times\sigma$ and $\phi^*(\om_Y\times\sigma)$ 
can be joined by a path of symplectic forms on $X\times S^2$.
\end{conjecture}

\noindent
The examples of Smith~\cite{SMITH1} in~\ref{ex:D} 
show that the $2$-sphere in Conjecture~\ref{con:don46}
cannot be replaced by the $2$-torus.  
The conjecture is nontrivial in either direction.  
When $X$ and $Y$ {\it are not} diffeomorphic but $X\times S^2$ 
and $Y\times S^2$ {\it are} diffeomorphic, the conjecture 
suggests that the two symplectic structures on these 
six-manifolds should still {\it remember} the 
differences in the smooth structures on $X$ and $Y$.
At the time of writing the only known methods for distinguishing
smooth structures on four-manifolds are the Donaldson invariants 
and the Seiberg--Witten invariants.  By Taubes--Seiberg--Witten 
theory~\cite{TAUBES3,TAUBES4,TAUBES5} the Seiberg--Witten
invariants can also be interpreted as symplectic invariants,
and it is conceivable that they give rise to a method 
for distinguishing the symplectic structures on 
the products with the two-sphere.
An example where $X$ and $Y$ {\it are not} diffeomorphic, 
$X\times S^2$ and $Y\times S^2$ {\it are} diffeomorphic, 
and the symplectic forms on $X\times S^2$ and $Y\times S^2$ 
can be distinguished by their Gromov--Witten 
invariants was found by Ruan~\cite{RUAN}
(see also~\cite[Example~9.7.1]{MS2}).
More examples along these lines were 
found by Ruan--Tian~\cite{RT} and Ionel--Parker~\cite{IP}.

\begin{conjecture}[{\bf Uniqueness conjecture}]
\label{con:torus}
Let $M$ be a closed hyperk\"ahler surface (i.e.\ a four-torus 
or a K3 surface) and let $a\in H^2(M;\R)$ be a cohomology 
class such that $a^2>0$. Then the space $\cS_a$ of symplectic 
forms on $M$ representing the class $a$ is connected.
\end{conjecture}

\noindent
The uniqueness problem for the four-torus is a longstanding 
open question in symplectic topology, which goes back at least
to the early eighties.  In~\cite{D4} Donaldson proposed a remarkable
geometric approach to the uniqueness question for symplectic forms 
on hyperk\"ahler surfaces which I explain next. 

\subsection*{The Donaldson geometric flow}

The starting point of Donaldson's approach is to view 
the space $\sF_a$ of diffeomorphisms from a symplectic 
four-manifold $(S,\sigma)$ to a hyperk\"ahler surface
$(M,\om_1,\om_2,\om_3,J_1,J_2,J_3)$, that pull back the cohomology 
class $a$ of $\om_1$ to that of $\sigma$, as an infinite dimensional
hyperk\"ahler manifold.  The group of symplectomorphisms
of $(S,\sigma)$ acts on $\sF_a$ by hyperk\"ahler isometries and the 
action of the subgroup of Hamiltonian symplectomorphisms is generated 
by a hyperk\"ahler moment map.  The negative gradient flow
of the square of the hyperk\"ahler moment map is
a parabolic type equation for a path of diffeomorphisms 
$f_t\in\sF_a$ depending on a real parameter $t$. It has the form
\begin{equation}\label{eq:sympflowf}
\p_tf_t = \sum_i J_idf_t\circ X_{H_i^{f_t}},\quad 
H_i^f := \frac{f^*\om_i\wedge\sigma}{\dvol_\sigma},\quad
\iota(X_H)\sigma = dH,
\end{equation}
where $\dvol_\sigma:=\tfrac12\sigma\wedge\sigma$.
Equation~\eqref{eq:sympflowf} is preserved 
by the symplectomorphism group of $(S,\sigma)$.
One can elmininate the action of the re\-pa\-ra\-me\-tri\-za\-tion 
group by pushing forward the sympectic form $\sigma$ under 
the diffeomorphisms~$f_t$. \linebreak  
The resulting path $\rho_t:=(f_t^{-1})^*\sigma\in\cS_a$  
satisfies the equation
\begin{equation}\label{eq:sympflow}
\p_t\rho_t = - \sum_i 
d\left(dK_i^{\rho_t}\circ J_i^{\rho_t}\right),\quad
K_i^\rho := \frac{\om_i\wedge\rho}{\dvol_\rho},\quad
\rho(J_i\cdot,\cdot) = \rho(\cdot,J_i^\rho\cdot).
\end{equation}
This is the {\bf Donaldson geometric flow}. It is the negative 
gradient flow of the energy functional $\cE:\cS_a\to\R$ defined by
\begin{equation}\label{eq:energy}
\cE(\rho) := \tfrac12\int_M\sum_i
\Abs{\frac{\om_i\wedge\rho}{\dvol_\rho}}^2\dvol_\rho
= 2\int_M \frac{\abs{\rho^+}^2}
{\abs{\rho^+}^2-\abs{\rho^-}^2}\dvol.
\end{equation}
The inner products on the tangent spaces
$T_\rho\cS_a=\left\{\widehat{\rho}\in\Om^2(M)\,|\,
\widehat{\rho}\mbox{ is exact}\right\}$
are associated to the norms
\begin{equation}\label{eq:inner}
\Norm{\widehat{\rho}}_\rho^2
:= \int_M\abs{X_{\widehat{\rho}}}^2\dvol_\rho,\qquad
\widehat{\rho}=-d\iota(X_{\widehat{\rho}})\rho,\qquad
\rho\wedge\inner{X_{\widehat{\rho}}}{\cdot}\in\im d.
\end{equation}
The gradient flow of the energy functional~\eqref{eq:energy} on $\cS_a$
with respect to the inner products determined by~\eqref{eq:inner}
is meaningful for any closed symplectic four-manifold $(M,\om)$, 
equipped with a background Riemannian metric.  It is convenient to 
choose the Riemannian metric on $M$ to be compatible with $\om$.
Then $\om$ is the unique global minimum of $\cE$ on $\cS_a$
and the Hessian of $\cE$ is positive definite at $\om$.    
In~\cite{D4} Donaldson also observed that the energy $\cE(\rho)$ 
controls the $L^1$-norm of $\rho$ and he proved, in the 
hyperk\"ahler case, that the Hessian at any other critical point
(if it exists) cannot be positive definite.  

\bigbreak

The study of equation~\eqref{eq:sympflow} poses deep 
and challenging analytical problems.  Already regularity 
and short time existence are nontrivial.  They are the subject 
of as yet unpublished work by Robin Krom~\cite{KROM},
which is based on ideas of Donaldson in~\cite{D5}.
The hope is that, in the hyperk\"ahler setting, one can establish 
long time existence and convergence for equation~\eqref{eq:sympflow}
and use this to prove that the space $\cS_a$ 
of symplectic forms in a fixed cohomology class is connected.
This hope is backed up by the fact that an analogous 
geometric flow approach in dimension two leads to the 
parabolic equation $\p_tu_t=d^*d u_t^{-1}$ (see~\cite{D4}).  
In this equation the proof of long time existence is based 
on the observation that the time derivative is nonpositive 
at each local maximum of~$u_t$ and is nonnegative 
at each local minimum of~$u_t$.

\medskip\noindent{\bf Acknowledgement.}
Thanks to Paul Biran, Simon Donaldson, Yael Kar\-shon,
Janko Lat\-schev, Dusa McDuff, and Stefano Vidussi 
for many helpful comments and suggestions.



\end{document}